\documentclass[a4paper,UKenglish,cleveref, autoref, thm-restate]{oasics-v2021}

\usepackage{mathtools}
\usepackage{tabularx}
\usepackage{booktabs}
\usepackage[short,c2,nocomma]{optidef}
\usepackage{etoolbox}
\usepackage{amsmath}
\usepackage{alphalph}
\usepackage{xspace}
\usepackage{pgfplots}
\pgfplotsset{compat=1.18}
\usepackage{standalone}
\usepackage{siunitx}

\allowdisplaybreaks

\AtBeginDocument{%
  \AtBeginEnvironment{subequations}{%
  }
}

\definecolor{PKdarkblue}{rgb}{0.121,0.47,0.705}
\definecolor{PKdarkred}{rgb}{0.89 0.102 0.109}
\definecolor{PKdarkgreen}{rgb}{0.2 0.627 0.172}
\definecolor{PKdarkorange}{rgb}{1 0.498 0}
\definecolor{PKdarkpurple}{rgb}{0.415 0.239 0.603}
\definecolor{PKlightgray}{rgb}{0.8 0.8 0.8}
\definecolor{PKdarkgray}{rgb}{0.5 0.5 0.5}
\definecolor{PKlightblue}{rgb}{0.651 0.807 0.89}
\definecolor{PKlightgreen}{rgb}{0.698 0.874 0.541}
\definecolor{PKlightorange}{rgb}{0.992 0.749 0.435}

\usepgfplotslibrary{statistics}
\pgfplotsset{compat=1.18,small,scale only axis,enlarge x limits={value=0.03,auto}}
\pgfplotscreateplotcyclelist{scatter-marks}{
  {solid,mark options={scale=.95,fill=PKdarkgreen},mark=*,PKdarkgreen},
  {solid,mark options={scale=.95,fill=PKdarkblue},mark=triangle*,PKdarkblue},
  {solid,mark options={scale=.95,fill=PKdarkorange},mark=diamond*,PKdarkorange},
}
\pgfplotscreateplotcyclelist{scatter-marks-2-cols}{
  {solid,mark options={scale=.95,fill=PKdarkgreen},mark=*,PKdarkgreen},
  {solid,mark options={scale=.95,fill=PKdarkgreen},mark=o,PKdarkgreen},
  {solid,mark options={scale=.95,fill=PKdarkblue},mark=triangle*,PKdarkblue},
  {solid,mark options={scale=.95,fill=PKdarkblue},mark=triangle,PKdarkblue},
  {solid,mark options={scale=.95,fill=PKdarkorange},mark=diamond*,PKdarkorange},
  {solid,mark options={scale=.95,fill=PKdarkorange},mark=diamond,PKdarkorange},
}
\pgfplotscreateplotcyclelist{scatter-2-marks}{
  {solid,mark options={scale=.95,fill=PKdarkpurple},mark=pentagon*,PKdarkpurple},
  {solid,mark options={scale=.95,fill=PKdarkorange},mark=diamond*,PKdarkorange},
}

\newcommand{\allK}{k \in K}
\newcommand{\allNR}[1]{#1 \in \NR}
\newcommand{\allP}[1]{#1 \in P}
\newcommand{\allH}[1]{#1 \in H}
\newcommand{\allHR}[1]{#1 \in \HR}
\newcommand{\allS}[1]{#1 \in S}
\newcommand{\allR}[1]{#1 \in R}
\newcommand{\subjH}[2]{\substack{j \in \HR:\\ (#1,#2) \in \ER}}
\newcommand{\vstart}[2]{\text{start}_{#1}^{#2}}
\newcommand{\vend}[2]{\text{end}_{#1}^{#2}}
\newcommand{\org}[1]{o_{{#1}}}
\newcommand{\dest}[1]{d_{{#1}}}
\newcommand{\pickup}[1]{\text{pickup}_{#1}}
\newcommand{\dropoff}[1]{\text{arrtime}_{#1}}
\newcommand{\maxridetime}[1]{L^{\max}_{#1}}
\newcommand{\precedenceSet}[1]{\text{Prec}_{#1}}

\newcommand{\lineproperty}{directionality property}
\newcommand{\distancereduction}{saved distance}
\newcommand{\baro}{\bar{o}}
\newcommand{\bard}{\bar{d}}
\newcommand{\barP}{\bar{P}}
\newcommand{\barD}{\bar{D}}
\newcommand{\tturn}{t_{\text{turn}}}
\newcommand{\dstart}{\delta_{\text{start}}}
\newcommand{\dend}{\delta_{\text{end}}}
\newcommand{\R}{R}
\newcommand{\GR}{G_{\R}}
\newcommand{\HR}{H_{\R}}
\newcommand{\ER}{E_{\R}}
\newcommand{\NR}{N_{\R}}
\newcommand{\Eturn}{E_{\text{turn}}}
\newcommand{\Qmax}{Q_{\max}}
\newcommand{\Sleft}{S^{\text{asc}}}
\newcommand{\Sright}{S^{\text{desc}}}
\newcommand{\depS}{\text{dep}}
\newcommand{\arrS}{\text{arr}}
\newcommand{\travelS}{x}
\newcommand{\stopS}{y}
\newcommand{\assignS}{\text{assign}}
\newcommand{\Rleft}{\R^{\text{asc}}}
\newcommand{\Rright}{\R^{\text{desc}}}
\newcommand{\Hleft}{H^{\text{asc}}}
\newcommand{\Hright}{H^{\text{desc}}}
\newcommand{\sizeS}{\sigma}
\newcommand{\sizeK}{\kappa}
\newcommand{\sizeR}{m}
\newcommand{\sizeH}{n}
\newcommand{\Tsubline}{T^+}

\pdfoutput=1
\hideOASIcs

\title{The Line-Based Dial-a-Ride Problem}

\author{Kendra Reiter\footnote{corresponding author}}{Department of Computer Science, University of Würzburg, Germany}{kendra.reiter@uni-wuerzburg.de}{https://orcid.org/0009-0004-7281-6516}{}
\author{Marie Schmidt}{Department of Computer Science, University of Würzburg, Germany}{marie.schmidt@uni-wuerzburg.de}{https://orcid.org/0000-0001-9563-9955}{}
\author{Michael Stiglmayr}{Department of Mathematics and Informatics, University of Wuppertal, Germany}{stiglmayr@uni-wuppertal.de}{https://orcid.org/0000-0003-0926-1584}{}

\authorrunning{K.~Reiter, M.~Schmidt, and M.~Stiglmayr}

\Copyright{Kendra Reiter, Marie Schmidt, and Michael Stiglmayr}

\ccsdesc[500]{Applied computing~Transportation}

\keywords{DARP, ridepooling, liDARP, public transport, on-demand}

\relatedversion{To appear in OASIcs, Volume 123, ATMOS 2024 (\url{https://www.dagstuhl.de/dagpub/978-3-95977-350-8}).}

\supplement{}
\supplementdetails[subcategory={Source Code \& Data}]{Software}{https://github.com/ReiterKM/liDARP}

\begin{document}

\maketitle

\begin{abstract}
On-demand ridepooling systems offer flexible services pooling multiple passengers into one vehicle, complementing traditional bus services. We propose a transportation system combining the spatial aspects of a fixed sequence of bus stops with the temporal flexibility of ridepooling. In the \emph{line-based Dial-a-Ride problem (liDARP)}, vehicles adhere to a fixed, ordered sequence of stops in their routes, with the possibility of taking shortcuts and turning if they are empty. We propose three MILP formulations for the liDARP with a multi-objective function balancing environmental aspects with customer satisfaction, comparing them on a real-world bus line. Our experiments show that the formulation based on an Event-Based graph is the fastest, solving instances with up to 50 requests in under one second. Compared to the classical DARP, the liDARP is computationally faster, with minimal increases in total distance driven and average ride times.
\end{abstract}

\section{Introduction}
\label{sec:intro}
Line-based bus services are able to pool a large number of transportation requests along popular trajectories, thus contributing towards reducing mobility-related emissions when people decide to take the bus instead of the car. However, when there is little demand (e.g., in rural areas or during off-peak periods), buses often run infrequently and almost empty, so that the described benefits do not materialize. Ridepooling approaches, where multiple passengers are pooled into a shared vehicle, accepting a slight detour compared to their direct route, are often proposed to complement line-based bus services for these scenarios. These approaches are inefficient where they do not succeed to pool requests sufficiently. 

We formalize a conceptual approach called \emph{line-based ridepooling} which combines the spatial aspect of a classical line-based bus service with the temporal flexibility of on-demand ridepooling, inspired by real-life examples, including the FLEX'HOP 72\footnote{\url{https://www.cts-strasbourg.eu/fr/se-deplacer/transport-a-la-demande/}} in France, and the NAHBUS\footnote{\url{https://www.nahbus.de/rufbus}} and the Rufbus\footnote{\url{https://rufbus.nordfriesland.de/Rufbus-N\"ordliches-u-s\"udliches-Nordfriesland/}} in Germany. In line-based ridepooling, we consider a fixed and ordered sequence of bus stops (defined, e.g., by a prior operating regular bus line) which we use as pick-up and drop-off locations. In contrast to a classical line-based bus, our vehicles have the flexibility to skip stops and take shortcuts, enabling tailored pick-up and drop-off times dependent on the specific customer requests. 
However, in contrast to ridepooling, we use the spatial structure of the given line as a sort of service promise, ensuring that passengers are only transported {towards their destination} along the line. In particular, vehicles may not turn with passengers on board. We call this the \emph{\lineproperty}.
We aim to achieve a transportation mode which is more efficient and provides a higher quality of service than the classical line-based bus, especially in areas with low demand or during off-peak times.

In this paper, we define a new optimization problem to 
serve passenger requests in line-based ridepooling. 
Due to its similarity to the general Dial-a-Ride problem (DARP), from which it differs by the \lineproperty\ and the underlying geography, we call this problem the \emph{line-based Dial-a-Ride problem (liDARP)}. 
Here, we study the static variant, where all requests are known ahead-of-time, and consider a homogeneous fleet of vehicles.

We observe that, due to the \lineproperty, each vehicle route can be decomposed into a number of sublines, separated by vehicle turns, with each transported passenger assigned to exactly one subline. We introduce and compare three mixed-integer linear programming (MILP) formulations for the liDARP that exploit this property. The first formulation, presented in \cref{sec:subline}, explicitly models sublines and the assignment of passengers to them. The second and third formulation, presented in \cref{sec:darp} and  \cref{sec:eb}, are based on Cordeau's classic 3-index Location-Based formulation \cite{cordeau_branch-and-cut_2006} and the Event-Based model introduced by Gaul et~al. in \cite{gaul_event-based_2022}. \Cref{sec:results} discusses computational results.

Our contribution is threefold: First, we present a general problem definition for the liDARP, an approach to organizing passenger transport with the potential to combine benefits from a line-based public transport and on-demand transportation. Second, we develop and present three MILP formulations
for the liDARP. Third, we compare these three models on synthetic test instances.

\section{Related Work}
Traditional modes for passenger transportation like the bus, metro, or train, operate based on lines (prescribing the sequence of stops visited) and timetables (prescribing the timing of each stop) or frequencies (prescribing the distance to be kept between individual vehicles on a line). 
While public transport planning often takes a network perspective, there are also many contributions that study timetabling or frequency setting on an individual line with the objective to find an optimal balance between service quality and operator cost, see, e.g., \cite{ibarra2015planning,liu2020robust} and the references therein. 
Gkiotsalitis et al.~\cite{gkiotsalitis_subline_2022} present a model that allows to establish regularly operating sublines within a longer line to deal with inhomogenous demand along the line.
Akta\c{s} et al.~\cite{aktas2022demand} study a situation where selected stops are assigned to an \emph{express service}, forming a shorter and quicker route. Their goal is to determine which vehicles should perform this express service during morning rush hour, based on expected demand.
While still a rather uncommon strategy during the \emph{planning} of public transport operations, short-turning and stop-skipping are common \emph{control}  strategies in transit systems to mitigate effects like vehicle bunching and overcrowding, see \cite{ibarra2015planning}.

The literature on Dial-a-Ride problems (synonymously called \emph{ridepooling}, \emph{on-demand} bus services, or \emph{demand-responsive transport}) is extensive, with in-depth overviews of the current state being provided by Cordeau and Laporte \cite{cordeau_Dial-a-Ride_2007} (until 2007) and Ho et~al. \cite{ho_survey_2018} (2007 until 2018).
Typography and variants are discussed in Molenbruch et~al. \cite{molenbruch_typology_2017}, whence this paper is concerned with the static, homogeneous, multi-objective approach, compromising the conflicting goals of system efficiency (including environmental aspects) and user experience.

Early approaches towards exact solution methods to the DARP were carried out by Psaraftis in \cite{psaraftis_dynamic_1980} and \cite{psaraftis_exact_1983}. Cordeau~\cite{cordeau_branch-and-cut_2006} proposes a 3-index arc-based mixed-binary linear program for the standard DARP, which was adapted to a 2-index formulation by R{\o}pke et~al.~\cite{ropke_models_2007}. R{\o}pke et~al.~propose a branch-and-cut approach which is tested on a large number of benchmark instances. 
Parragh~\cite{parragh_introducing_2011} constructs further valid inequalities related to capacity restrictions, integrating these into a branch-and-cut framework as well as a variable neighborhood search heuristic, based on both the 3-index and 2-index formulations.
Gschwindt and Irnich~\cite{gschwind2015effective} develop an exact branch-and-cut-and-price approach, which solves all instances of the benchmark set introduced by \cite{ropke_models_2007} exactly. Recently, Rist and Forbes~\cite{rist_new_2021} propose a branch-and-cut framework where a DARP route is broken in multiple \emph{fragments}, which are paths between a request's pick-up and drop-off where the vehicle has a non-empty load. Then, a route is created as a combination of fragments. Gaul et~al.~\cite{gaul_event-based_2022} propose a new MILP~formulation, relying on an \emph{Event-Based} graph with nodes representing pick-up/drop-off events denoting a feasible user allocation of the corresponding vehicle and edges connecting feasible transitions between events.

Next to exact methods, many papers consider heuristic solution methods to solve the DARP, including metaheuristics such as simulated annealing \cite{baugh_intractability_1998, reinhardt_synchronized_2013}, adaptive large neighborhood search \cite{parragh_hybrid_2013, pfeiffer_alns_2022,ropke_adaptive_2006}, and tabu search \cite{attanasio_parallel_2004, cordeau_tabu_2003}. 

We are aware of only three publications where the DARP is studied in combination with an underlying line structure:
Archetti et~al.~\cite{archetti_complexity_2011} restate and prove results from the dissertation of Busch \cite{busch_vehicle_1991}, showing that the Vehicle Routing Problem on the line is NP-hard, both with an unlimited and a limited fleet of fixed capacity. A complexity classification of DARP variants has been proposed by de~Paepe et~al.~\cite{de2004computer}, establishing a scheme akin to scheduling problems. They examine variants on the line geography, showing that the DARP on the line with one vehicle of capacity one is solvable in polynomial time. The DARP on a line with multiple homogeneous vehicles of fixed capacity $\geq 1$ is NP-complete, which has been shown by Bjelde et~al.~\cite{bjelde2020tight} based on a reduction from the \textsc{Circular Arc Coloring} problem. All three papers are focused on exploring the complexity of the problem, where they consider only the special case with equally spaced stations and do not allow for shortcuts.

\section{Problem Description}\label{sec:prob_descr}
We consider a set of $\sizeK$ vehicles of capacity $\Qmax$ that operate on a bus line, specified by a sequence of bus stops $H= (1, \ldots, n)$, to transport $\sizeR$ stop-to-stop passenger requests $\R$. 

The vehicles do not need to traverse the whole line in each route, but are allowed to take short-cuts (including skipping stops at which no passenger wants to board or alight), to wait, and to turn at any stop, the latter of which may not be done with passengers on board. In this way, we guarantee that the \emph{\lineproperty} is fulfilled, i.e., each passenger, at all times, is transported towards their direction with respect to the sequence of stops defined by the bus line.
Pairwise (time) distances $t_{i,j}$ between all stops $i,j\in H$ are given, with $t_{i,i} := \tturn$ denoting the turn time at $i \in H$. These distances respect the triangle inequality.

Each request $r\in \R$ specifies an origin stop $o_r \in H$, a destination stop $d_r \in H$, a time window $[e_r,\ l_r]$, a load (number of passengers in the request) $q_r$, and a service time $b_r$ for boarding and alighting.
We assume that boarding is synchronous, i.e., if one request's destination is another request's origin, and both requests are transported by the same vehicle, we require that the first request alights before the second boards. This reflects the widely accepted standard boarding procedure on public transit systems. Furthermore, passengers do not transfer between vehicles.

We make two \emph{service promises} to our accepted passengers regarding 1)~their total travel time and 2)~their waiting time that have to be respected. For the former, we guarantee that the passenger's total ride time $L_r$ will not exceed the time needed to travel the direct route (between their origin and destination) by a pre-specified factor, the \emph{excess factor} $\alpha$, i.\,e., $L_r^{\text{max}} := \alpha \cdot t_{o_r, d_r}$. For the latter, we ensure that the actual pick-up (resp.\ drop-off) time is not more than $\beta$ minutes later (resp.\ earlier) than the specified earliest pick-up (resp.\ latest drop-off) time.

Our objective is to create a reliable service
for customers and, at the same time, integrate environmental aspects by reducing emissions compared to passengers travelling in their own vehicles. Therefore, our objective function is composed of two weighted components: the number of \emph{accepted} passengers and the \emph{\distancereduction}\ (i.\,e., the difference between the sum of direct distances between all origins and destinations and the total distance driven by our vehicles), which we want to maximize. The optimization problem now consists of deciding which passenger requests are accepted, and which are rejected, to assign accepted requests to one of the $\sizeK$ vehicles, and to plan the routes of these vehicles.

The above-defined problem is a variant of the Dial-a-Ride problem: removing the restriction that vehicles may only turn without passengers on board, it reduces to a (standard) DARP.
Given that our vehicle's operations are constrained by the line, we call our problem the line-based Dial-a-Ride problem (liDARP).

\section{MILP Formulations for the liDARP}
The underlying line structure, which defines an order of bus stops, combined with the \lineproperty, allows us to divide the route of each vehicle into a number of \emph{sublines}: a sequence of stops at which a vehicle stops to pick-up or drop-off passengers. The first subline is initialized when a vehicle starts its route and a new subline starts after each of the vehicle's turns. We split the set of sublines $S$ into \emph{ascending sublines} ($\Sleft$), travelling from a stop $i$ to $j$ with $i < j$, and \emph{descending sublines} ($\Sright$), travelling in the opposite direction.

Similarly, we divide the passenger requests $r \in \R$ into \emph{ascending requests} ($\R^{\text{asc}}$) and \emph{descending requests} ($\R^{\text{desc}}$). As passengers may not be on board when the vehicle turns, each accepted request can be assigned to exactly one subline, with ascending requests assigned to ascending sublines and descending requests assigned to descending sublines.

In \cref{sec:subline}, we exploit these properties to propose a \emph{Subline-Based} MILP for the liDARP, explicitly modelling sublines and passenger assignments to sublines. In \cref{sec:darp} and \cref{sec:eb}, we show that the sublines can also be used to simplify MILP formulations for the standard DARP. 

Note that the sublines in the liDARP are similar to the so-called \emph{fragments} proposed by Rist and Forbes~\cite{rist_new_2021} in their branch-and-cut approach for the DARP. Namely, a subline can be further subdivided into fragments, which start with a pick-up node and end when the vehicle is empty.

\subsection{Subline-Based Formulation}\label{sec:subline}
The \emph{Subline-Based formulation} relies on the concept of a subline. Given the set of vehicles $K$, we assign each $k \in K$ a set of sublines $S$ and use binary variables $y_i^{s,k}$ to indicate whether subline $s$ of vehicle $k$ stops at bus stop $i$. The route of every subline $s$ of vehicle $k$ is encoded by the binary variables $\travelS_{i,j}^{s,k}$ that denote the path between bus stops $i, j \in H$. Depending on the direction of the subline $s$, these only need to be defined for $i \leq j$ or $j \leq i$, respectively. Sublines are computed on a vehicle-basis and symmetry breaking constraints are defined on the vehicle's index to remove alternative solutions with equal objectives.

Flow conservation constraints ensure these routes are consistent in each subline and between consecutive sublines. We track the start and end stop of each subline with binary variables $\travelS_{i,i}^{s,k}$, which correspond to sublines $s$ of vehicle $k$ turning at bus stop $i$. 

Requests are assigned to sublines using binary variables $\assignS_r ^{s,k}$ to indicate if request $r$ is transported by subline $s$ of vehicle $k$. Note that these only need to be created for pairs $(r,s)$ with both $r$ and $s$ travelling in the same direction. We ensure each passenger is picked up at most once, with the corresponding subline stopping at both the origin and destination stop. The underlying line structure determines each subline's pick-up and drop-off sequence, allowing capacity constraints to be expressed solely in variables $\assignS_r ^{s,k}$.

Continuous variables $\arrS_i^{s,k}$ and $\depS_i^{s,k}$ model the arrival and departure time of subline $s$ of vehicle $k$ at bus stop $i$, respectively. We introduce constraints to ensure the stopping time is sufficiently long for all assigned passengers who are boarding or alighting the vehicle at a station to do so. Similarly, we ensure that the time between departure at a bus stop $i$ and arrival at the next bus stop $j$ on the vehicle's route is equal to $t_{i,j}$. The departure and arrival times are further constrained by the fact that, when a request is assigned to a subline, the corresponding departure and arrival times have to respect the request's time window and associated service promises. For this, we track the pick-up and arrival time of every request $r$.

An overview of parameters and variables, and the full model are given in \cref{appendix:milp:subline}.

\subsection{Location-Based Formulation}\label{sec:darp}
Inspired by mathematical programming formulation for the traveling salesperson and vehicle routing problems, Cordeau~\cite{cordeau_branch-and-cut_2006} models the Dial-a-Ride problem using a graph where nodes represent origin $o_r$ and destination $d_r$ locations of requests $r$ and arcs represent direct connections between locations. In principle, traveling between any pair of locations is possible, though many arcs can be removed in a pre-processing step based on time constraints.

To account for the \lineproperty, we can modify Cordeau's DARP formulation for use in the liDARP as follows: we treat the requests as (general) DARP input, 
adding modifications to respect the line precedence and to prevent that vehicles turn with passengers on board. Observing that a vehicle may only turn after drop-off location or before a pick-up location, we introduce additional nodes at these bus stops which are used to start or end a turn, allowing us to model our problem based on fewer arcs than the general DARP. 
A \emph{start-turn} node denotes that a vehicle is turning at a pick-up node (and then starting a new subline), while an \emph{end-turn} node denotes that a vehicle is turning after a drop-off node (and ending the current subline). Similar to the general DARP, many arcs can be excluded by pre-processing based on time windows and service constraints.

Binary variables $x_{i,j}^k$ model whether vehicle $k$ travels from node $i$ to node $j$ with flow constraints ensuring feasible operations. Additional variables and big-$M$-constraints are needed to keep track of the time at which locations are visited and of vehicle load, so that the requirements with respect to time windows, service promises, and vehicle capacity are ensured. Technical details and the full model are given in \cref{appendix:milp:darp}, where we use strengthening techniques based on \cite{desrochers_improvements_1991}.

\subsection{Event-Based Formulation}\label{sec:eb}
Encoding feasible user allocations in vehicles as nodes, and constructing only edges between feasible connections, Gaul et al.~\cite{gaul_event-based_2022} propose the \emph{Event-Based} graph as a basis to formulate the general DARP as a MILP. Every node in the Event-Based graph represents a $\Qmax$-tuple, where the first entry represents the most recent action: a pick-up $r^+$ or drop-off $r^-$ of a request $r$, and the remaining entries denote the other passengers on board. The node $\mathbf{0}$ is used to denote the depot.
Finding feasible vehicle routes for the DARP can then be interpreted as a minimum cost circulation flow problem with the additional constraints that each passenger cannot be picked up more than once and that time windows and service constraints need to be respected.

While the number of events grows exponential with $\Qmax$, many nodes can already be excluded during the construction of the Event-Based graph due to incompatibility of time windows and service constraints. The \lineproperty\ allows us to further reduce this set. In particular, requests $i$ and $j$ cannot be part of the same event if
\begin{itemize}
    \item one of them is ascending and the other is descending, or
    \item both requests are ascending (resp. descending) and the request with a later (resp. earlier) starting station cannot board a vehicle with the other already on board due to time window constraints. 
\end{itemize}
The Event-Based model for the DARP can be directly applied to the liDARP by an adapted construction of the Event-Based graph, distinguishing ascending and descending events and connecting only events that preserve the \lineproperty, hence we do not re-state the formulation here. Moreover, Gaul et~al.~\cite{gaul_tight_2023} proposed further arc eliminations.

\section{Computational Experiments}\label{sec:results}
In this section, we present numerical experiments for the liDARP on synthetic benchmark instances. For all experiments, we set the passenger load $q_r = 1$ and service time $b_r = \SI{3}{\minute}$ for all $r \in \R$. The service promise parameters were set to a maximum waiting time of $\beta = \SI{15}{\minute}$ and a maximum exceedance of direct ride time by $\alpha = 3$. The objective function weights were chosen to be $c_1 = 10$ for the number of accepted passenger and $c_2 = 1$ for the \distancereduction\ for the computational results.

The models for the Subline-Based and Location-Based formulation were implemented in Python 3.11 using Gurobi 10.0. The Event-Based formulation was implemented in C++ 17 using CPLEX 22.1, based on the code by Gaul et~al.~\cite{gaul_solving_2021}.
The computations are carried out using a 12th Gen Intel Core i7-1260P CPU, running at \SI{2.10}{\giga\hertz} with \SI{32}{\giga\byte} RAM. For all runs, we set the solver timeout to \SI{60}{\minute} and repeated the calculation five times, averaging the runtimes.

\subsection{Benchmark Instances}
\begin{figure}[htb!]
    \centering
    \includegraphics[width=.7\linewidth, clip]{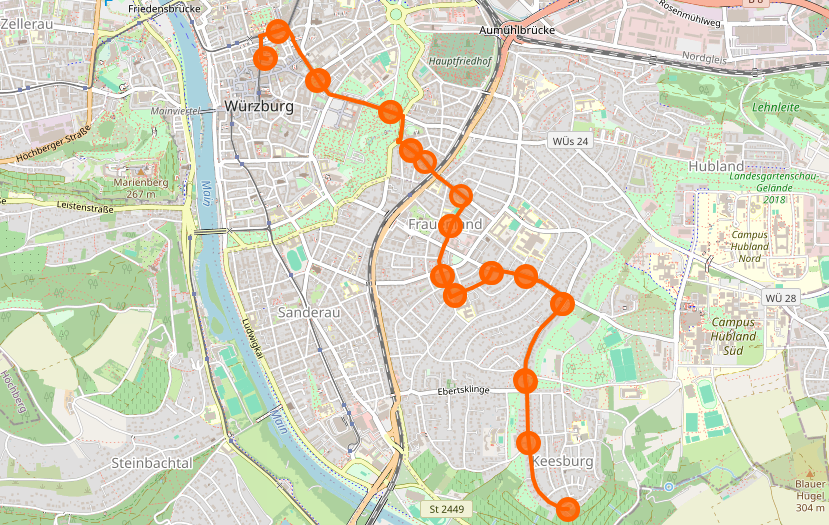}
    \caption{Route of bus line 6 in Würzburg, Germany, from \cite{OpenStreetMap}.}
    \label{fig:route_line_6}
\end{figure}

We create new benchmark instances specifically for the liDARP using the existing bus stops of bus line 6 in Würzburg, Germany, as pictures in \cref{fig:route_line_6}, with 16 stops connecting the city center to a residential area. We calculate bus stop distances using OpenStreetMap~\cite{OpenStreetMap}, assuming vehicles can take shortcuts and rounding to the nearest minute.

For the given bus stops $H$, we generate requests that uniformly choose a pick-up and drop-off stop in $H$. We generate an equal amount of requests per time window type, picking the (earliest) pick-up time or (latest) drop-off time uniformly in the interval $[0,480]$, corresponding to an operation of 8 hours. The vehicles have a capacity of $\Qmax \in \{3,6\}$ and take $\tturn = \SI{3}{\minute}$ to turn.

We generated 14 instances on the given sequence $H$, varying from 16 requests with 2 vehicles to 50 requests with 5 vehicles, following the sizing of the well-known benchmark instances by Cordeau~\cite{cordeau_branch-and-cut_2006} for the classical DARP.
The instance names consist of a prefix `w' (for Würzburg), followed by two numbers, where the first indicates the number of vehicles and the second denotes the number of requests.

\subsection{Results}
In this section, we first compare the three proposed MILP formulations using the benchmark bus line test instances. Second, we assess the trade-off between environmental savings and customer satisfaction. Lastly, we compare the liDARP model to the classical DARP model to evaluate its competitiveness.

\subsubsection*{Computational Time}
\begin{figure}[hbt!]
    \centering
    \begin{tikzpicture}
        \begin{semilogyaxis}[
            table/col sep=comma,
            width=10cm,height=5cm,
            symbolic x coords={w2-16, w2-20,w2-24,w3-18,w3-24,w3-30,w3-36,w4-16,w4-24,w4-32,w4-40,w4-58,w5-40,w5-50},
            xticklabel style={rotate=45,yshift=0.6em,xshift=-0.4em},
            cycle list name=scatter-marks-2-cols,
            legend columns=2,
            legend style={at={(0.03,0.65)},anchor=west,nodes={scale=0.8, transform shape}},
            legend cell align=left,
            xlabel={instance},
            ylabel={runtime [s]}
            ]
            \addplot+[only marks] table[x=instance,y=subline-based] {average_runtime.dat};
            \addlegendentry{}
            \addplot+[only marks] table[x=instance,y=event-based] {SB_Q6_avg_runtime.dat};
            \addlegendentry{Subline-Based}
            \addplot+[only marks] table[x=instance,y=location-based] {average_runtime.dat};
            \addlegendentry{}
            \addplot+[only marks] table[x=instance,y=event-based] {LB_Q6_avg_runtime.dat};
            \addlegendentry{Location-Based}
            \addplot+[only marks] table[x=instance,y=event-based] {average_runtime.dat};
            \addlegendentry{}
            \addplot+[only marks] table[x=instance,y=event-based] {EB_Q6_avg_runtime.dat};
            \addlegendentry{Event-Based}
            \draw[thick,PKlightgray,dashed] (rel axis cs:0,0.916) -- (rel axis cs:1,0.916);
        \end{semilogyaxis}
    \end{tikzpicture}
    \caption{Runtime of the three formulations on the benchmark instances. Solid markers denote $\Qmax=3$, unfilled markers $\Qmax=6$, and the dashed line marks the solver timeout.}
    \label{fig:runtime}
\end{figure}
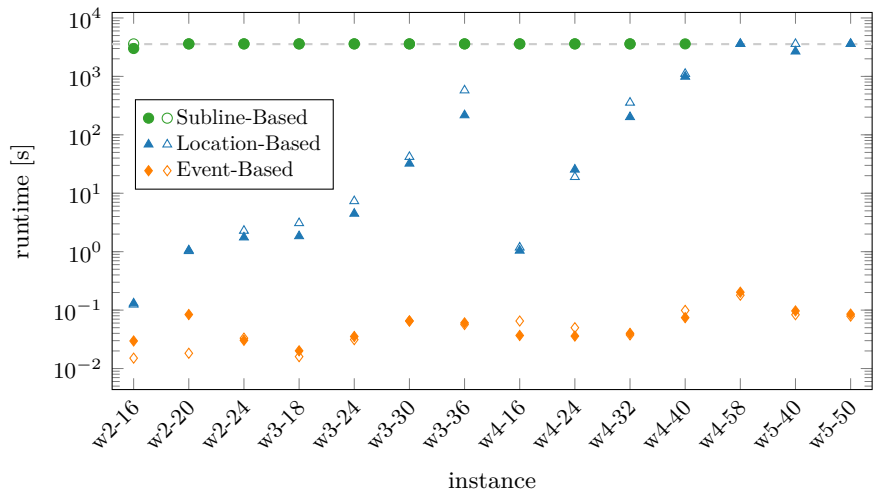

s\cref{fig:runtime} shows the computational time per benchmark instance for all three formulations. 

The results clearly show that the Event-Based model outperforms the Location-Based and Subline-Based models in all instances, for either capacity. For $\Qmax = 3$, the Subline-Based model reached the timeout for all but one instance, w2-16, and we were not able to compute a solution for the three largest instances due a lack of available memory. For $\Qmax = 6$, the largest possible instance the model could solve was w4-32.

Observing the Location-Based model, we see a step-structure, where the runtime significantly increases with the number of requests, then decreases as the instances switch to the next-largest number of vehicles. We see a similar effect in the Event-Based model between instances w4-58 and w5-40. This is supported by the model size differences reported in \cref{tab:results:num_const_and_vars}, which strongly correlate to the number of passenger requests.

In the two instances where the Location-Based model reached timeout for $\Qmax=3$, namely w4-58 and w5-50, we note that although the achieved a relative MIP gap at timeout was greater than 1, the objective value found was within \qty{8}{\percent} and even \qty{0}{\percent} of the optimum, respectively.

\begin{table}[hbt!]
    \centering
    \caption{Number of constraints and variables for the benchmark test instances with $\Qmax = 3$, ordered by the number of requests. SB = Subline-Based, LB = Location-Based, EB = Event-Based.}\label{tab:results:num_const_and_vars}
    \begin{tabular}{lrrrrrrrrr}
        & \multicolumn{3}{c}{Num. Constraints} & \multicolumn{3}{c}{Num. Boolean Var.} & \multicolumn{3}{c}{Num. Cont. Var.} \\
        \cmidrule(lr){2-4}\cmidrule(lr){5-7}\cmidrule(lr){8-10}
        Inst. & SB & LB & EB & SB & LB & EB & SB & LB & EB  \\
        \midrule
        w2-16 & \num{45042} & \num{553} & \num{317} & \num{20030} & \num{968} & \num{189} & \num{1952} & \num{150} & \num{53} \\
        w4-16 & \num{90388} & \num{787} & \num{335} & \num{40188} & \num{1868} & \num{196} & \num{3872} & \num{154} & \num{56} \\
        w3-18 & \num{83315} & \num{752} & \num{401} & \num{35841} & \num{1731} & \num{238} & \num{3276} & \num{170} & \num{64} \\
        w2-20 & \num{67290} & \num{689} & \num{497} & \num{28062} & \num{1372} & \num{294} & \num{2440} & \num{186} & \num{76} \\
        w2-24 & \num{96434} & \num{825} & \num{746} & \num{38462} & \num{2062} & \num{430} & \num{2928} & \num{222} & \num{106} \\
        w3-24 & \num{144443} & \num{998} & \num{830} & \num{57117} & \num{3168} & \num{449} & \num{4368} & \num{224} & \num{119} \\
        w4-24 & \num{197732} & \num{1171} & \num{769} & \num{79228} & \num{4116} & \num{437} & \num{5808} & \num{226} & \num{109} \\
        w3-30 & \num{227963} & \num{1244} & \num{956} & \num{87213} & \num{4473} & \num{610} & \num{5460} & \num{278} & \num{122} \\
         w4-32 & \num{352756} & \num{1555} & \num{1174} & \num{133500} & \num{6792} & \num{706} & \num{7744} & \num{298} & \num{149} \\
        w3-36 & \num{340895} & \num{1490} & \num{1627} & \num{127749} & \num{6564} & \num{935} & \num{6552} & \num{332} & \num{199} \\
        w4-40 & \num{577636} & \num{1939} & \num{1867} & \num{214204} & \num{10668} & \num{1117} & \num{9680} & \num{370} & \num{223} \\
        w5-40 & \num{722789} & \num{2224} & \num{1657} & \num{267755} & \num{12945} & \num{1065} & \num{12080} & \num{372} & \num{188} \\
        w5-50 & \num{1227649} & \num{2774} & \num{2419} & \num{453155} & \num{20885} & \num{1602} & \num{15100} & \num{462} & \num{249} \\
        w4-58 & \num{1400268} & \num{2803} & \num{3282} & \num{517020} & \num{21976} & \num{2154} & \num{14036} & \num{532} & \num{326} \\
        \bottomrule
        \end{tabular}
\end{table}

The Subline-Based model requires a significantly higher number of resources, with the number of constraints and boolean variables exceeding those of both the Location-Based and Event-Based models by factors of 100 and 10, respectively. Notably, the Location-Based model uses 10 times more boolean variables than the Event-Based model, but both require a similar amount of continuous variables.

In general, this result is not surprising, as Gaul et~al.~\cite{gaul_event-based_2022} demonstrated the computational efficiency of the Event-Based graph in a MILP formulation for the classical DARP, as many complicating constraints are implicitly encoded in the underlying network structure.

All following experiments are carried out with $\Qmax=3$.

\subsubsection*{Trade-off Analysis}
\begin{figure}[htb!]
    \centering
    \begin{tikzpicture}
        \begin{axis}[
            table/col sep=comma,
            width=10cm,
            height=5cm,
            legend style={at={(0.04,0.95)},anchor=north west,nodes={scale=0.8, transform shape}},
            legend cell align=left,
            xlabel={share of accepted passengers},
            ylabel={\distancereduction},
            xtick={0,10,20,30,40,50,60,70,80,90,100},
            xmin = 28,
            xticklabel={$\pgfmathprintnumber{\tick}\%$}
            ]
            \addlegendimage{only marks, mark=*,PKdarkgreen},
            \addlegendentry{environmental focus},
            \addlegendimage{mark=*,PKdarkblue,only marks},
            \addlegendentry{equal weights},
            \addlegendimage{mark=*,PKdarkorange,only marks},
            \addlegendentry{customer focus},
            \addplot+[only marks,
            scatter/classes={
            0={mark=*,fill={PKdarkgreen},draw={PKdarkgreen}},
            1={mark=triangle*,fill={PKdarkgreen},draw={PKdarkgreen}},
            2={mark=square*,fill={PKdarkgreen},draw={PKdarkgreen}},
            3={mark=diamond*,fill={PKdarkgreen},draw={PKdarkgreen}},
            4={mark=pentagon*,fill={PKdarkgreen},draw={PKdarkgreen}},
            5={mark=o,fill={PKdarkgreen},draw={PKdarkgreen}},
            6={mark=triangle,fill={PKdarkgreen},draw={PKdarkgreen}},
            7={mark=square,fill={PKdarkgreen},draw={PKdarkgreen}},
            8={mark=diamond,fill={PKdarkgreen},draw={PKdarkgreen}},
            9={mark=pentagon,fill={PKdarkgreen},draw={PKdarkgreen}},
            10={mark=star,fill={PKdarkgreen},draw={PKdarkgreen}},
            11={mark=+,fill={PKdarkgreen},draw={PKdarkgreen}},
            12={mark=x,fill={PKdarkgreen},draw={PKdarkgreen}},
            13={mark=asterisk,fill={PKdarkgreen},draw={PKdarkgreen}}
            },
            scatter src=explicit symbolic,
            scatter
            ] table[x=Ecological-emphasis_x,y=Ecological-emphasis_y,meta=perc_accepted] {vary_weights_perc_accepted.dat};
            \addplot+[only marks,
            scatter/classes={
            0={mark=*,fill={PKdarkblue},draw={PKdarkblue}},
            1={mark=triangle*,fill={PKdarkblue},draw={PKdarkblue}},
            2={mark=square*,fill={PKdarkblue},draw={PKdarkblue}},
            3={mark=diamond*,fill={PKdarkblue},draw={PKdarkblue}},
            4={mark=pentagon*,fill={PKdarkblue},draw={PKdarkblue}},
            5={mark=o,fill={PKdarkblue},draw={PKdarkblue}},
            6={mark=triangle,fill={PKdarkblue},draw={PKdarkblue}},
            7={mark=square,fill={PKdarkblue},draw={PKdarkblue}},
            8={mark=diamond,fill={PKdarkblue},draw={PKdarkblue}},
            9={mark=pentagon,fill={PKdarkblue},draw={PKdarkblue}},
            10={mark=star,fill={PKdarkblue},draw={PKdarkblue}},
            11={mark=+,fill={PKdarkblue},draw={PKdarkblue}},
            12={mark=x,fill={PKdarkblue},draw={PKdarkblue}},
            13={mark=asterisk,fill={PKdarkblue},draw={PKdarkblue}}
            },
            scatter src=explicit symbolic,
            scatter
            ] table[x=Equal weights_x,y=Equal weights_y,meta=perc_accepted] {vary_weights_perc_accepted.dat};
            \addplot+[only marks,
            scatter/classes={
            0={mark=*,fill={PKdarkorange},draw={PKdarkorange}},
            1={mark=triangle*,fill={PKdarkorange},draw={PKdarkorange}},
            2={mark=square*,fill={PKdarkorange},draw={PKdarkorange}},
            3={mark=diamond*,fill={PKdarkorange},draw={PKdarkorange}},
            4={mark=pentagon*,fill={PKdarkorange},draw={PKdarkorange}},
            5={mark=o,fill={PKdarkorange},draw={PKdarkorange}},
            6={mark=triangle,fill={PKdarkorange},draw={PKdarkorange}},
            7={mark=square,fill={PKdarkorange},draw={PKdarkorange}},
            8={mark=diamond,fill={PKdarkorange},draw={PKdarkorange}},
            9={mark=pentagon,fill={PKdarkorange},draw={PKdarkorange}},
            10={mark=star,fill={PKdarkorange},draw={PKdarkorange}},
            11={mark=+,fill={PKdarkorange},draw={PKdarkorange}},
            12={mark=x,fill={PKdarkorange},draw={PKdarkorange}},
            13={mark=asterisk,fill={PKdarkorange},draw={PKdarkorange}}
            },
            scatter src=explicit symbolic,
            scatter
            ] table[x=Customer-emphasis_x,y=Customer-emphasis_y, meta=perc_accepted] {vary_weights_perc_accepted.dat};
            \draw[thick,PKlightgray,dashed] (0,0) -- (100,0);
            \draw[thin, PKdarkgray] (40.0,7.0) -- (75.0, 0.00) -- (100.0,-22.0);
        \end{axis}
    \end{tikzpicture}
    \caption{Objective Function parameters for varying weights. Each marker shape represents one benchmark instance.}
    \label{fig:vary-weights-accepted}
\end{figure}
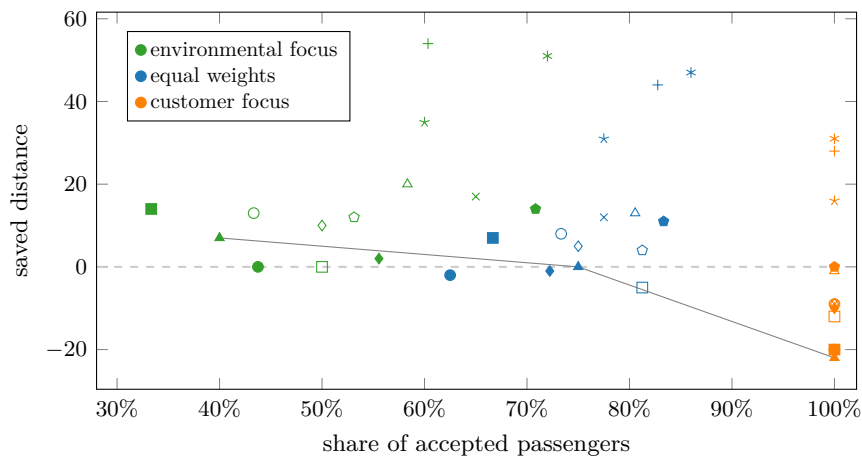

To evaluate the trade-off between environmental savings and customer attractiveness in our chosen objective function, we compare three different settings: in the \emph{environmentally focused} setting, we use objective weights $c_1 = 1$, $c_2 = 10$, placing an emphasis on the distance saved, whilst in the \emph{customer focused} setting, we use weights $c_1 = 10, c_2 = 1$, emphasising the number of transported passengers. We also include a setting with equal weights as a base case. The trade-off between objective function components is visualized in \cref{fig:vary-weights-accepted}. 

Note that both objective functions are considered as maximizing objectives. A positive \distancereduction\ is preferred, as this corresponds to more direct passenger kilometers saved than total routing costs accumulated. All instances in \cref{fig:vary-weights-accepted} were solved to optimality.

We observe that, in the customer-focused setting, all passengers are accepted in all instances. In the environmental-focused setting, the \distancereduction\ is always non-negative. We have highlighted the approximated Pareto front for a specific instance, w2-20, represented by triangles, by connecting the markers corresponding to the three obtained weighted-sum solutions in \cref{fig:vary-weights-accepted} to better illustrate the trade-offs. The saved distance decreases from \qty{7}{\minute} to \qty{-22}{\minute} between the environmental-focused and the customer-focused setting, while the share of accepted passengers increases from \qty{40}{\percent} to \qty{100}{\percent}. The same pattern can be observed for all other benchmark instances. Hence, the proposed objective function is capable of capturing multiple needs and can be adjusted accordingly, dependent on the chosen application.

\subsubsection*{DARP versus liDARP}
\begin{figure}[htb!]
    \centering
    \begin{tikzpicture}
        \begin{semilogyaxis}[
            table/col sep=comma,
            width=10cm,
            height=5cm,
            symbolic x coords={w2-16, w2-20,w2-24,w3-18,w3-24,w3-30,w3-36,w4-16,w4-24,w4-32,w4-40,w4-58,w5-40,w5-50,w5-60,w6-48,w6-60,w6-72,w7-56,w7-70,w7-84,w8-64,w8-80,w8-96,w9-72,w9-90,w9-108,w10-80,w10-100,w10-120,w11-88,w11-110,w11-132},
            xticklabel style={rotate=45,yshift=0.6em,xshift=-0.4em},
            cycle list name=scatter-2-marks,
            legend style={at={(0.03,0.65)},anchor=west,nodes={scale=0.8, transform shape}},
            legend cell align=left,
            xlabel={instance},
            ylabel={runtime [s]},
            ymax = 15000,
            minor x tick num = 1,
            ]
            \addplot+[only marks] table[x=instance,y=event-based] {darp_avg_runtime.dat};
            \addlegendentry{DARP};
            \addplot+[only marks] table[x=instance,y=event-based] {lidarp_0tturn_avg_runtime.dat}; 
            \addlegendentry{liDARP}
            \draw[thick,PKlightgray,dashed] (rel axis cs:0,0.915) -- (rel axis cs:1,0.915);
        \end{semilogyaxis}
    \end{tikzpicture}
    \caption{Runtime of the DARP and liDARP on the extended benchmark instances. The dashed line marks the solver timeout.}
    \label{fig:lidar-vs-darp}
\end{figure}
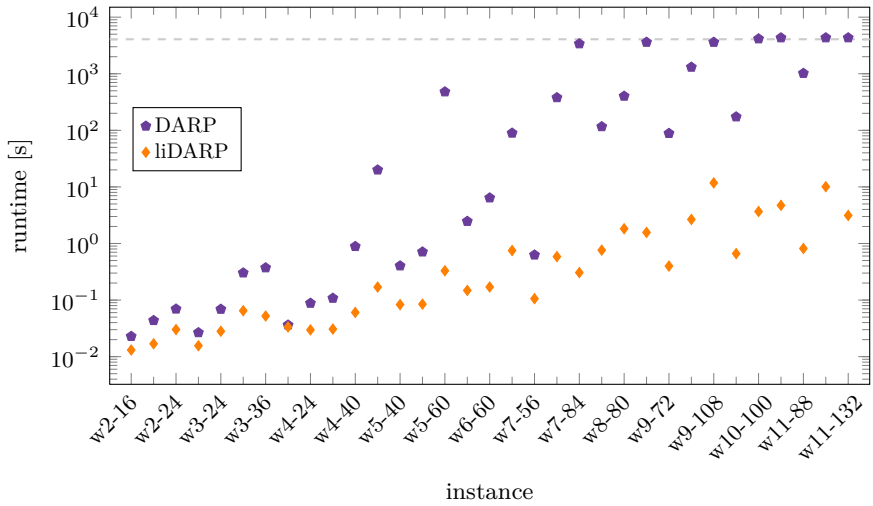

Lastly, we compare the liDARP formulation to the classical DARP formulation, where vehicles are allowed to take any route between passenger without needing to adhere to the line structure. We use the above-introduced benchmark instances, extending these to cases with up to 11 vehicles and 132 requests (similar to the extended benchmark set introduced by R{\o}pke et~al.~\cite{ropke_models_2007}). We set $\tturn = \SI{0}{\minute}$ for both formulations and set the solver timeout to \qty{60}{\minute}. The computational time for both models is shown in \cref{fig:lidar-vs-darp}, averaged over five runs. The Event-Based model was used to produce the liDARP results.

\begin{figure}[htb!]
    \centering
    \begin{tikzpicture}
        \begin{axis}[
            table/col sep=comma,
            width=10cm,
            height=5cm,
            symbolic x coords={w2-16, w2-20,w2-24,w3-18,w3-24,w3-30,w3-36,w4-16,w4-24,w4-32,w4-40,w4-58,w5-40,w5-50,w5-60,w6-48,w6-60,w6-72,w7-56,w7-70,w7-84,w8-64,w8-80,w8-96,w9-72,w9-90,w9-108,w10-80,w10-100,w10-120,w11-88,w11-110,w11-132},
            xticklabel style={rotate=45,yshift=0.6em,xshift=-0.4em},
            cycle list name=scatter-2-marks,
            legend style={at={(0.03,0.65)},anchor=north west,nodes={scale=0.8, transform shape}},
            legend cell align=left,
            xlabel={instance},
            ylabel={saved distance},
            ytick = {0,50,100},
            ymin=-50,
            ymax=150,
            minor x tick num = 1,
            ]
            \addplot+[only marks] table[x=instance,y=event-based] {darp_avg_saved_distance.dat};
            \addlegendentry{DARP};
            \addplot+[only marks] table[x=instance,y=event-based] {lidarp_0tturn_avg_saved_distance.dat}; 
            \addlegendentry{liDARP}
            \draw[thick,PKlightgray,dashed] (rel axis cs:0,0.25) -- (rel axis cs:1,0.25);
        \end{axis}
    \end{tikzpicture}
    \caption{Saved distance of the liDARP and DARP on the extended benchmark instances.}
    \label{fig:lidar-vs-darp:saved_distance}
\end{figure}
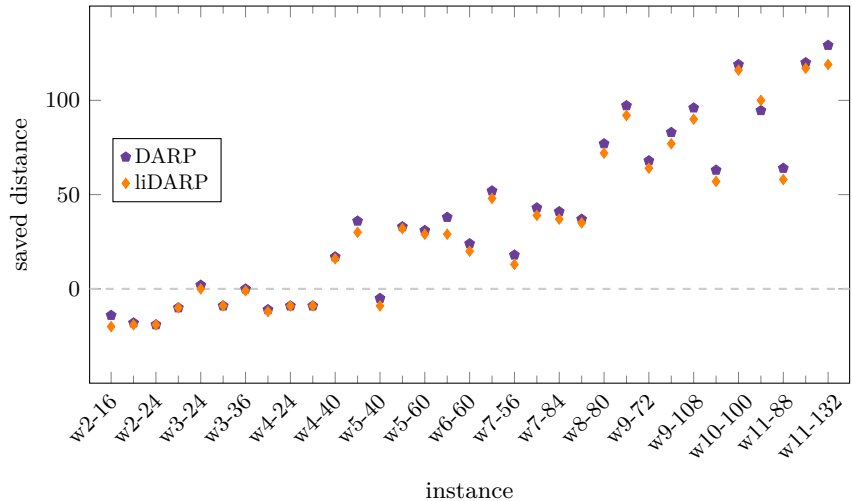

We observe that the liDARP model is faster in all instances. While both models' computational time increases with the number of requests, the liDARP was able to solve even the largest instances with over 100 requests in less than \SI{10}{\second}, while the DARP model was aborted at timeout. Examining the objective values, both models accepted all requests in all instances, while they differ in the saved distance, which is visualized in \cref{fig:lidar-vs-darp:saved_distance}. The DARP achieved a higher \distancereduction\ in all but one instance, w10-100, which was aborted at timeout, with the average deviation being \qty{3}{\minute} and the maximum deviation being a saving of \qty{10}{\minute} in instance w11-132. Both models use all available vehicles for all instances.

The average ride time, measured as the time between pick-up and drop-off of each request, was marginally higher in the DARP model, with an average increase of \qty{0.14}{\minute}. The difference in the average share of empty mileage, defined as the fraction of empty mileage over the total distance, is less than $0.02$ on average. Similarly, the difference in average detour, which is the fraction of passenger distance driven to the shortest distances between pick-up and drop-off, is less than $0.025$ on average. 

Examining each DARP solution, we count the number of requests which are, at least for a portion of their trip, travelling \emph{away} from their destination, i.e., violating the \lineproperty. The only instances where there are no such violations are w3-18, w3-30, and w4-32. On average, \qty{7.7}{\percent} of passengers travel in the opposite direction for at least a portion of their trip, with the largest amount being \qty{18.9}{\percent} in instance w9-90. We hypothesize that these routes will likely be viewed as unnecessary by customers, even if they are the most efficient amongst all possible connections.

While there is a difference in saved distance and ride time, passengers do not have to accept significant detours when using the liDARP compared to the DARP model, and the vehicles travel without passengers for a similar amount of time.
\section{Conclusion}

We present the line-based Dial-a-Ride problem (liDARP), wherein ridepooling vehicles operate on-demand on a sequence of bus stations, adhering to the \lineproperty, time and capacity constraints, and our service promises. The efficiency of this approach is validated through numerical computations on benchmark instances derived from on a real-life bus line in Würzburg, Germany. Our multi-objective approach successfully balances environmental concerns, by reducing the total distance travelled compared to individual passenger trips, and the attractiveness, measured in the total amount of passengers accepted. 

The advantages of the liDARP additionally include the possibility to use existing infrastructure, being based on a sequence of bus stops, which may increase utilization in off-peak times or in areas with low demand. Thus, we provide a system that can serve as an alternative to public transport, providing flexibility to its customers to improve attractiveness.

Our approach can be extended to inhomogeneous vehicles with varying capacities, to determine where each vehicle is best allocated to serve the customer base. Future research could explore varying demand scenarios, such as operating a feeder line to a train station or considering rush-hour induced fluctuations in demand. Finally, this paper focuses on the static variant of the liDARP, where all requests are known in advance. Investigating the liDARP under unknown and dynamic demand may provide insights into its practical applicability and competitiveness in a real-world setting.

\bibliography{lidarp-article}

\begin{thebibliography}{10}

\bibitem{aktas2022demand}
Dilay Aktas, Pieter Vansteenwegen, and Kenneth S{\"o}rensen.
\newblock A demand-responsive bus system for peak hours with capacitated
  vehicles.
\newblock In {\em Proc.\ 11th Triennial Symposium on Transportation Analysis
  conference (TRISTAN XI)}, Mauritius Island, 2022. TRISTAN.

\bibitem{archetti_complexity_2011}
Claudia Archetti, Dominique Feillet, Michel Gendreau, and M.~Grazia~Speranza.
\newblock Complexity of the {VRP} and {SDVRP}.
\newblock {\em Transportation Research Part C: Emerging Technologies},
  19(5):741--750, August 2011.
\newblock \href {https://doi.org/10.1016/j.trc.2009.12.006}
  {\path{doi:10.1016/j.trc.2009.12.006}}.

\bibitem{attanasio_parallel_2004}
Andrea Attanasio, Jean-François Cordeau, Gianpaolo Ghiani, and Gilbert
  Laporte.
\newblock Parallel {Tabu} search heuristics for the dynamic multi-vehicle
  dial-a-ride problem.
\newblock {\em Parallel Computing}, 30(3):377--387, March 2004.
\newblock \href {https://doi.org/10.1016/j.parco.2003.12.001}
  {\path{doi:10.1016/j.parco.2003.12.001}}.

\bibitem{baugh_intractability_1998}
John~W. Baugh, Gopala Krishna~Reddy Kakivaza, and John~R. Stone.
\newblock Intractability of the {Dial}-a-{Ride} {Problem} and a
  {Multiobjective} {Solution} {Using} {Simulated} {Annealing}.
\newblock {\em Engineering Optimization}, 30(2):91--123, February 1998.
\newblock \href {https://doi.org/10.1080/03052159808941240}
  {\path{doi:10.1080/03052159808941240}}.

\bibitem{bjelde2020tight}
Antje Bjelde, Jan Hackfeld, Yann Disser, Christoph Hansknecht, Maarten Lipmann,
  Julie Mei{\ss}ner, Miriam Schl{\"o}ter, Kevin Schewior, and Leen Stougie.
\newblock Tight bounds for online tsp on the line.
\newblock {\em ACM Transactions on Algorithms (TALG)}, 17(1):1--58, 2021.
\newblock \href {https://doi.org/10.1145/3422362} {\path{doi:10.1145/3422362}}.

\bibitem{busch_vehicle_1991}
Ingrid Busch.
\newblock {\em Vehicle routing on acyclic networks}.
\newblock Dissertation, The Johns Hopkins University, Baltimore, Maryland,
  1991.

\bibitem{cordeau_branch-and-cut_2006}
Jean-Fran{\c{c}}ois Cordeau.
\newblock A {Branch}-and-{Cut} {Algorithm} for the {Dial}-a-{Ride} {Problem}.
\newblock {\em Operations Research}, 54(3):573--586, 2006.
\newblock \href {https://doi.org/10.1287/opre.1060.0283}
  {\path{doi:10.1287/opre.1060.0283}}.

\bibitem{cordeau_tabu_2003}
Jean-Fran{\c{c}}ois Cordeau and Gilbert Laporte.
\newblock A tabu search heuristic for the static multi-vehicle dial-a-ride
  problem.
\newblock {\em Transportation Research Part B: Methodological}, 37(6):579--594,
  2003.
\newblock \href {https://doi.org/10.1016/S0191-2615(02)00045-0}
  {\path{doi:10.1016/S0191-2615(02)00045-0}}.

\bibitem{cordeau_Dial-a-Ride_2007}
Jean-Fran{\c{c}}ois Cordeau and Gilbert Laporte.
\newblock The dial-a-ride problem: models and algorithms.
\newblock {\em Annals of Operations Research}, 153(1):29--46, 2007.
\newblock \href {https://doi.org/10.1007/s10479-007-0170-8}
  {\path{doi:10.1007/s10479-007-0170-8}}.

\bibitem{de2004computer}
Willem~E. de~Paepe, Jan~Karel Lenstra, Jiri Sgall, Ren\'{e}~A. Sitters, and
  Leen Stougie.
\newblock Computer-aided complexity classification of dial-a-ride problems.
\newblock {\em INFORMS Journal on Computing}, 16(2):120--132, 2004.
\newblock \href {https://doi.org/10.1287/ijoc.1030.0052}
  {\path{doi:10.1287/ijoc.1030.0052}}.

\bibitem{desrochers_improvements_1991}
Martin Desrochers and Gilbert Laporte.
\newblock Improvements and extensions to the {Miller}-{Tucker}-{Zemlin} subtour
  elimination constraints.
\newblock {\em Operations Research Letters}, 10(1):27--36, 1991.
\newblock \href {https://doi.org/10.1016/0167-6377(91)90083-2}
  {\path{doi:10.1016/0167-6377(91)90083-2}}.

\bibitem{gaul_tight_2023}
Daniela Gaul, Kathrin Klamroth, Christian Pfeiffer, Arne Schulz, and Michael
  Stiglmayr.
\newblock A {Tight} {Formulation} for the {Dial}-a-{Ride} {Problem}, 2023.
\newblock \href {https://arxiv.org/abs/2308.11285} {\path{arXiv:2308.11285}}.

\bibitem{gaul_solving_2021}
Daniela Gaul, Kathrin Klamroth, and Michael Stiglmayr.
\newblock {Solving the Dynamic Dial-a-Ride Problem Using a Rolling-Horizon
  Event-Based Graph}.
\newblock In {\em 21st Symposium on Algorithmic Approaches for Transportation
  Modelling, Optimization, and Systems (ATMOS 2021)}, volume~96 of {\em Open
  Access Series in Informatics (OASIcs)}, pages 8:1--8:16, Dagstuhl, Germany,
  2021. Schloss Dagstuhl -- Leibniz-Zentrum f{\"u}r Informatik.
\newblock \href {https://doi.org/10.4230/OASIcs.ATMOS.2021.8}
  {\path{doi:10.4230/OASIcs.ATMOS.2021.8}}.

\bibitem{gaul_event-based_2022}
Daniela Gaul, Kathrin Klamroth, and Michael Stiglmayr.
\newblock Event-based {MILP} models for ridepooling applications.
\newblock {\em European Journal of Operational Research}, 301(3):1048--1063,
  2022.
\newblock \href {https://doi.org/10.1016/j.ejor.2021.11.053}
  {\path{doi:10.1016/j.ejor.2021.11.053}}.

\bibitem{gkiotsalitis_subline_2022}
Konstantinos Gkiotsalitis, Marie Schmidt, and Evelien van~der Hurk.
\newblock Subline frequency setting for autonomous minibusses under demand
  uncertainty.
\newblock {\em Transportation Research Part C: Emerging Technologies},
  135:103492, 2022.
\newblock \href {https://doi.org/10.1016/j.trc.2021.103492}
  {\path{doi:10.1016/j.trc.2021.103492}}.

\bibitem{gschwind2015effective}
Timo Gschwind and Stefan Irnich.
\newblock Effective handling of dynamic time windows and its application to
  solving the dial-a-ride problem.
\newblock {\em Transportation Science}, 49(2):335--354, 2015.
\newblock \href {https://doi.org/10.1287/trsc.2014.0531}
  {\path{doi:10.1287/trsc.2014.0531}}.

\bibitem{ho_survey_2018}
Sin~C. Ho, Wai~Yuen Szeto, Yong-Hong Kuo, Janny~M.Y. Leung, Matthew Petering,
  and Terence~W.H. Tou.
\newblock A survey of dial-a-ride problems: {Literature} review and recent
  developments.
\newblock {\em Transportation Research Part B: Methodological}, 111:395--421,
  2018.
\newblock \href {https://doi.org/10.1016/j.trb.2018.02.001}
  {\path{doi:10.1016/j.trb.2018.02.001}}.

\bibitem{ibarra2015planning}
Omar~J. Ibarra-Rojas, Felipe Delgado, Ricardo Giesen, and Juan~Carlos
  Mu{\~n}oz.
\newblock Planning, operation, and control of bus transport systems: A
  literature review.
\newblock {\em Transportation Research Part B: Methodological}, 77:38--75,
  2015.
\newblock \href {https://doi.org/10.1016/j.trb.2015.03.002}
  {\path{doi:10.1016/j.trb.2015.03.002}}.

\bibitem{liu2020robust}
Pei Liu, Marie Schmidt, Qingxia Kong, Joris~Camiel Wagenaar, Lixing Yang, Ziyou
  Gao, and Housheng Zhou.
\newblock A robust and energy-efficient train timetable for the subway system.
\newblock {\em Transportation Research Part C: Emerging Technologies},
  121:102822, 2020.
\newblock \href {https://doi.org/10.1016/j.trc.2020.102822}
  {\path{doi:10.1016/j.trc.2020.102822}}.

\bibitem{molenbruch_typology_2017}
Yves Molenbruch, Kris Braekers, and An~Caris.
\newblock Typology and literature review for dial-a-ride problems.
\newblock {\em Annals of Operations Research}, 259(1):295--325, 2017.
\newblock \href {https://doi.org/10.1007/s10479-017-2525-0}
  {\path{doi:10.1007/s10479-017-2525-0}}.

\bibitem{OpenStreetMap}
{OpenStreetMap contributors}.
\newblock {Planet dump retrieved from https://planet.osm.org }.
\newblock \url{ https://www.openstreetmap.org }, 2017.

\bibitem{parragh_introducing_2011}
Sophie~N. Parragh.
\newblock Introducing heterogeneous users and vehicles into models and
  algorithms for the dial-a-ride problem.
\newblock {\em Transportation Research Part C: Emerging Technologies},
  19(5):912--930, 2011.
\newblock \href {https://doi.org/10.1016/j.trc.2010.06.002}
  {\path{doi:10.1016/j.trc.2010.06.002}}.

\bibitem{parragh_hybrid_2013}
Sophie~N. Parragh and Verena Schmid.
\newblock Hybrid column generation and large neighborhood search for the
  dial-a-ride problem.
\newblock {\em Computers \& Operations Research}, 40(1):490--497, 2013.
\newblock \href {https://doi.org/10.1016/j.cor.2012.08.004}
  {\path{doi:10.1016/j.cor.2012.08.004}}.

\bibitem{pfeiffer_alns_2022}
Christian Pfeiffer and Arne Schulz.
\newblock An {ALNS} algorithm for the static dial-a-ride problem with ride and
  waiting time minimization.
\newblock {\em OR Spectrum}, 44(1):87--119, 2022.
\newblock \href {https://doi.org/10.1007/s00291-021-00656-7}
  {\path{doi:10.1007/s00291-021-00656-7}}.

\bibitem{psaraftis_dynamic_1980}
Harilaos~N. Psaraftis.
\newblock A {Dynamic} {Programming} {Solution} to the {Single} {Vehicle}
  {Many}-to-{Many} {Immediate} {Request} {Dial}-a-{Ride} {Problem}.
\newblock {\em Transportation Science}, 14(2):130--154, 1980.
\newblock \href {https://doi.org/10.1287/trsc.14.2.130}
  {\path{doi:10.1287/trsc.14.2.130}}.

\bibitem{psaraftis_exact_1983}
Harilaos~N. Psaraftis.
\newblock An {Exact} {Algorithm} for the {Single} {Vehicle} {Many}-to-{Many}
  {Dial}-{A}-{Ride} {Problem} with {Time} {Windows}.
\newblock {\em Transportation Science}, 17(3):351--357, 1983.
\newblock \href {https://doi.org/10.1287/trsc.17.3.351}
  {\path{doi:10.1287/trsc.17.3.351}}.

\bibitem{reinhardt_synchronized_2013}
Line~Blander Reinhardt, Tommy Clausen, and David Pisinger.
\newblock Synchronized dial-a-ride transportation of disabled passengers at
  airports.
\newblock {\em European Journal of Operational Research}, 225(1):106--117,
  2013.
\newblock \href {https://doi.org/10.1016/j.ejor.2012.09.008}
  {\path{doi:10.1016/j.ejor.2012.09.008}}.

\bibitem{rist_new_2021}
Yannik Rist and Michael~A. Forbes.
\newblock A {New} {Formulation} for the {Dial}-a-{Ride} {Problem}.
\newblock {\em Transportation Science}, 55(5):1113--1135, 2021.
\newblock \href {https://doi.org/10.1287/trsc.2021.1044}
  {\path{doi:10.1287/trsc.2021.1044}}.

\bibitem{ropke_models_2007}
Stefan Ropke, Jean-François Cordeau, and Gilbert Laporte.
\newblock Models and branch-and-cut algorithms for pickup and delivery problems
  with time windows.
\newblock {\em Networks}, 49(4):258--272, 2007.
\newblock \href {https://doi.org/10.1002/net.20177}
  {\path{doi:10.1002/net.20177}}.

\bibitem{ropke_adaptive_2006}
Stefan Ropke and David Pisinger.
\newblock An {Adaptive} {Large} {Neighborhood} {Search} {Heuristic} for the
  {Pickup} and {Delivery} {Problem} with {Time} {Windows}.
\newblock {\em Transportation Science}, 40(4):455--472, 2006.
\newblock \href {https://doi.org/10.1287/trsc.1050.0135}
  {\path{doi:10.1287/trsc.1050.0135}}.

\end{thebibliography}
\appendix
\section{Variable Overview}

\begin{table}[htb!]
    \caption{Summary of liDARP parameters.}
    \label{tab:summary:all:params}
    \begin{tabularx}{\textwidth}{ll}
        \toprule
        Notation & Definition \\
        \midrule
        $H$ & set of bus stops, $\{1, \ldots, \sizeH\}$ \\
        $K$ & set of vehicles, $\{1, \ldots, \sizeK\}$ \\
        $\R$ & set of passenger requests, $\{1, \ldots, \sizeR\}$  \\
        $\Rleft$ & set of passenger requests travelling in ascending direction \\
        $\Rright$ & set of passenger requests travelling in descending direction \\
        $o_r$ & origin stop of passenger request $r$ \\
        $d_r$ & destination stop of passenger request $r$ \\
        $q_r$ & load of request $r$ \\
        $e_r$ & earliest departure time of request $r$ \\
        $l_r$ & latest arrival time of request $r$ \\
        $b_r$ & service time for request $r$ \\
        $t_{i,j}$ & travel time from bus stop $i$ to $j$ \\
        $\alpha$ & service promise constant relating to maximum ride time, excess factor\\
        $\beta$ & service promise constant relating to maximum wait time\\
        $L_r$ & total ride time of request $r$\\
        $\maxridetime{r}$ & maximum ride time of request $r$, dependent on $\alpha$ and $\beta$ \\
        $\tturn$ & time it takes for a vehicle to turn around  \\
        $c_1, c_2$ & objective weights \\
        \bottomrule
    \end{tabularx}
\end{table}

\section{MILP Formulations}

\subsection{Subline-Based Formulation}\label{appendix:milp:subline}

In this section, we present the MILP model for the Subline-Based formulation introduced in \cref{sec:subline}. All parameters and variables are summarized in \cref{tab:summary:sublines:params}.

\begin{table}[hbt]
    \centering
    \caption{Summary of notation for the Subline-Based model.}
    \label{tab:summary:sublines:params}
    \begin{tabularx}{\textwidth}{ll}
    \toprule
    Notation & Definition \\
    \midrule
    \multicolumn{2}{l}{Parameters} \\
    \midrule
    $S$ & set of sublines, $\{1, \ldots, \sizeS\}$ \\
    $\Sleft$ & set of sublines travelling in ascending direction \\
    $\Sright$ & set of sublines travelling in descending direction \\
    \midrule
    \multicolumn{2}{l}{Binary Decision Variables} \\
    \midrule
    $\assignS_r^{s,k}$ & $1$ if request $r$ is assigned to subline $s$ of vehicle $k$\\
    $\stopS_{i}^{s,k}$  & $1$ if subline $s$ of vehicle $k$ visits node $i$\\
    $\vstart{i}{k}$ & $1$ if node $i$ is the start node of vehicle $k$\\
    $\vend{i}{k}$ & $1$ if node $i$ is the end node of vehicle $k$ \\
    $\travelS_{i,j}^{s,k}$ & $1$ if node $j$ is visited immediately after node $i$ on subline $s$ of vehicle $k$\\
    $\travelS_{i,i}^{s,k}$ & $1$ if vehicle $k$ turns at node $i$ after executing subline $s$\\
    $w_{r_1, r_2}$ & $1$ if requests $r_1$ and $r_2$ are on the same subline of the same vehicle \\
    $z_k$ & $1$ if vehicle $k$ is in use\\
    \midrule
    \multicolumn{2}{l}{Continuous Decision Variables} \\
    \midrule
    $\depS_i^{s,k}$ & departure time of subline $s$ of vehicle $k$ at node $i$\\
    $\arrS_i^{s,k}$ & arrival time of subline $s$ of vehicle $k$ at node $i$ \\
    $\pickup{r}$ & pick-up time of request $r$\\
    $\dropoff{r}$ & drop-off time of request $r$\\
    \bottomrule
    \end{tabularx}
\end{table}

In this model, we explicitly model the path of every subline of each vehicle, using binary variables $\stopS_{i}^{s,k}$ to denote if subline $s$ of vehicle $k$ stops at bus station $i$ and binary variables $\travelS_{i,j}^{s,k}$ to denote if the direct path from station $i$ to station $j$ is used. Further binary variables $\vstart{i}{k}$ and $\vend{i}{k}$ denote if vehicle $k$ starts and ends at station $i$, respectively. Then, by tracking the turning stations of every subline $s$, i.e., stations $i$ where $\travelS_{i,i}^{s,k} = 1$, we track the start and end station of every subline. For every variable which references both passengers and sublines, such as $\assignS_r^{s,k}$, we only create those variables where the passenger and subline are travelling in the same direction.

To ensure our model respects the boarding precedence (requests which are alighting leave the vehicle before those boarding can enter), we construct the following three subsets for every request $i \in \R$:
\begin{itemize}
    \item $\precedenceSet{i}^{o,o}$, containing all requests $j$ which have the same origin as $i$, are travelling in the same direction, and should board before $i$,
    \item $\precedenceSet{i}^{o,d}$, containing all requests $j$ whose destination is at $i$'s origin, are travelling in the same direction, and should alight before $i$ boards, and
    \item $\precedenceSet{i}^{d,d}$, containing all requests $j$ which have the same destination as $i$, are travelling in the same direction, and should alight before $i$.
\end{itemize}

Note that the set $\precedenceSet{i}^{d,o}$, which denotes all origin requests that need to be served before request $i$ is dropped-off, is not created as it is always empty due to the boarding assumption.

Furthermore, we track which passengers are pooled together on the same subline of the same vehicle with binary variables $w_{r_1, r_2}$. These are required, together with the precedence sets, to ensure we allow for sufficient boarding and alighting times per passenger at every stop. To model these, we use big-M constraints with $M_1 := \max_{r \in \R} \maxridetime{r}$.

Binary variables $z_k$ denote that vehicle $k$ is used by our solution. Here, we use big-M constraints with $M_2 := \max_{r \in R}l_r + (2 \cdot \sizeS - 1) \cdot b_r$, where $\sizeS$ denotes the number of sublines,  to ensure the arrival and departure time variables for each vehicle are only set if they are also used.

To strengthen the model, we enforce that a vehicle's end stop is placed after it has turned twice, i.e., after two consecutive sublines start and end at the same stop. Then, every following subline is empty and turns at the same stop. This reduces the number of possible solutions with the same objective value.

We define $T^+ := \max_{r \in R} \maxridetime{r} + (\sizeS - 1) \tturn$ to be the end of service, i.e., when all vehicles end their operation at the latest.

We note that the model's size is dependent on the choice of $\sizeS$, the number of sublines, which is hard to choose. We set $\sizeS = 2 \cdot \sizeR$ for all experiments presented here.

The full Subline-Based model is given by:
\begin{subequations}
    \begin{align}
        \max_x\  &c_1 \cdot \Biggl(\sum_{s \in S}\sum_{k \in K}\sum_{r \in \R} \assignS_r^{s,k} \cdot t_{\org{r}, \dest{r}}  - \sum_{k \in K}\sum_{s \in S}\sum_{\substack{(i,j) \in H \times H: \\ i \neq j}} \travelS_{i,j}^{s,k} \cdot t_{i,j}  \Biggr)  \notag \\
        & + c_2 \cdot \sum_{s \in S}\sum_{k \in K}\sum_{r \in \R} \assignS_r^{s,k} \notag \\
        \text{s.t.}\ & \sum_{i \in H} \vstart{i}{k} \leq 1 \qquad\forall\ \allK \label{model:subline:subline:1}\\ 
        &{\sum_{i \in H} \vend{i}{k}}{ = \sum_{i \in H} \vstart{i}{k}}{\qquad\forall\ \allK\label{model:subline:subline:2}}\\ 
        &{\sum_{i \in H} \travelS_{ii}^{s,k}}{= \sum_{i \in H}\vstart{i}{k}}{\qquad\forall\ \allK, s \in S \setminus \{\sizeS\}\label{model:subline:subline:3}}\\ 
        &{\travelS_{i,i}^{\sizeS, k}}{= 0}{\qquad\forall\ \allH{i}, \allK \label{model:subline:subline:4}}\\ 
        &{\stopS_i^{1,k}}{= \vstart{i}{k} + \sum_{j < i}  \travelS_{j,i}^{1,k}}{\qquad\forall\ \allH{i}, \allK \label{model:subline:subline:5}}\\ 
        &{\stopS_i^{s,k}}{= \travelS_{i,i}^{s-1, k} + \sum_{j < i} \travelS_{j,i}^{s,k}}{\qquad\forall\ \allH{i}, \allK, s \in \Sleft \setminus \{1\} \label{model:subline:subline:6}}\\ 
        &{\stopS_i^{s,k}}{= \travelS_{i,i}^{s-1, k} + \sum_{j > i} \travelS_{j,i}^{s,k}}{\qquad\forall\ \allH{i}, \allK, s \in \Sright}\label{model:subline:subline:7}\\ 
        &{\stopS_i^{s,k}}{= \travelS_{i,i}^{s,k} + \sum_{j > i}\travelS_{i,j}^{s,k}}{\qquad\forall\ \allH{i}, \allK, s \in \Sleft}\label{model:subline:subline:8}\\ 
        &{\stopS_i^{s,k}}{= \travelS_{i,i}^{s,k} + \sum_{j < i}\travelS_{i,j}^{s,k}}{\qquad\forall\ \allH{i}, \allK, s \in \Sright \setminus \{\sizeS\}}\label{model:subline:subline:9}\\ 
        &{\stopS_i^{\sizeS, k}}{= \vend{i}{k} + \sum_{j < i}\travelS_{i,j}^{\sizeS, k}}{\qquad\forall\ \allH{i}, \allK}\label{model:subline:subline:10}\\ 
        &{z_k}{ = \sum_{i \in H} \vstart{i}{k}}{\qquad\forall\ \allK}\label{model:subline:subline:11}\\ 
        &{\sum_{r : \org{r} \leq i, \dest{r} > i} \assignS_r^{s,k}}{\leq \Qmax}{\qquad\forall\  \allK, s \in \Sleft, i = 1, \ldots, n-1}\label{model:subline:subline:12}\\ 
        &{\sum_{r : \org{r} \geq i, \dest{r} < i} \assignS_r^{s,k}}{ \leq \Qmax}{\qquad\forall\ \allK, s \in \Sright, i = 2, \ldots, n}\label{model:subline:subline:13}\\
        &{z_k}{\geq \stopS_i^{s,k}}{\qquad\forall\ \allH{i}, \allS{s}, \allK}\label{model:subline:subline:14}\\ 
        &{ z_k}{\geq \travelS_{i,j}^{s,k}}{\qquad \forall (i,j) \in E, \allS{s}, \allK}\label{model:subline:subline:15}\\
        &{z_k}{\geq \assignS_r^{s,k}}{\qquad\forall\ \allR{r}, \allS{s}, \allK}\label{model:subline:subline:16}\\
        &{M_2 \cdot z_k}{\geq \arrS_i^{s,k}}{\qquad\forall\ \allH{i}, \allS{s}, \allK}\label{model:subline:subline:17}\\ 
        &{M_2 \cdot z_k}{\geq \depS_i^{s,k}}{ \qquad\forall\ \allH{i}, \allS{s}, \allK}\label{model:subline:subline:18}\\
        &{\vstart{i}{k}}{\geq \sum_{j \in H: j < i} \vstart{j}{k+1}}{\qquad\forall\ \allH{i}, k \in \{1, \ldots, \sizeK - 1\}}\label{model:subline:subline:19}\\
        &{\depS_i^{s,k}}{ \geq \arrS_i^{s,k}}{\qquad\forall\ \allH{i}, \allS{s}, \allK}\label{model:subline:time:1}\\ 
        &{\arrS_j^{s,k}}{ \geq \depS_i^{s,k} + t_{i,j} \cdot \travelS_{i,j}^{s,k}}{\qquad\forall (i,j) \in H \times H \text{ with } i < j, s \in \Sleft, \allK}\label{model:subline:time:2}\\ 
        &{\arrS_j^{s,k}}{\geq \depS_i^{s,k} + t_{i,j} \cdot \travelS_{i,j}^{s,k}}{\qquad\forall (i,j) \in H \times H \text{ with } i > j, s \in \Sright, \allK}\label{model:subline:time:3}\\
        &{\arrS_i^{s,k}}{\geq \depS_i^{s-1, k} + \tturn \cdot \travelS_{i,i}^{s-1, k}}{\qquad\forall\ \allH{i}, s \in \Sleft \setminus \{1\}, \allK}\label{model:subline:time:4}\\ 
        &{\arrS_i^{s,k}}{ \geq \depS_i^{s-1, k} + \tturn \cdot \travelS_{i,i}^{s-1, k}}{\qquad\forall\ \allH{i}, s \in \Sright, \allK}\label{model:subline:time:5}\\
        &{ \arrS_i^{s,k}}{\geq 0}{\qquad\forall\ \allH{i}, \allS{s}, \allK}\label{model:subline:time:6}\\ 
        &{\sum_{k \in K} \sum_{s \in S} \assignS_r^{s,k}}{\leq 1}{\qquad\forall\ \allR{r}}\label{model:subline:pax:subline:1}\\ 
        &{2 \cdot \assignS_r^{s,k}}{\leq \stopS_{\org{r}}^{s,k} + \stopS_{\dest{r}}^{s,k}}{\qquad\forall\ \allK, \allS{s}, \allR{r}}\label{model:subline:pax:subline:2}\\ 
        &{\assignS_r^{s,k} \cdot (e_r + b_{o_r})}{\leq \depS_{\org{r}}^{s,k}}{\qquad\forall\ \allK, \allS{s}, \allR{r}}\label{model:subline:pax:subline:5}\\ 
        &{l_r + \Tsubline \cdot (1 - \assignS_r^{s,k})}{\geq \arrS_{\dest{r}}^{s,k}}{\qquad\forall\ \allR{r}, \allS{s}, \allK}\label{model:subline:pax:subline:6}\\ 
        &{w_{r_i, r_j}}{\leq  \assignS_{r_i}^{s,k}}{\qquad\forall\ \allS{s}, \allK, \allR{r_i, r_j}}\label{model:subline:pax:subline:7}\\ 
        &{\assignS_{r_i}^{s,k} + \assignS_{r_j}^{s,k} - 1}{\leq w_{r_i, r_j}}{\qquad\forall\ \allS{s}, \allK, \allR{r_i, r_j}}\label{model:subline:pax:subline:8}\\
        &{\pickup{r}}{\geq \sum_{k \in K} \sum_{s \in S} \assignS_r^{s,k} \cdot \arrS_{\org{r}}^{s,k}}{\qquad\forall\ \allK, \allS{s}, \allR{r}}\label{model:subline:pax:time:1}\\ 
        &{\dropoff{r}}{\geq \sum_{k \in K} \sum_{s \in S} \assignS_r^{s,k} \cdot \arrS_{\dest{r}}^{s,k}}{\qquad\forall\ \allK, \allS{s}, \allR{r}}\label{model:subline:pax:time:2}\\ 
        &{b_{r_i} + \pickup{r_i} -\pickup{r_j}}{\leq M_1 \cdot \Biggl( 1 - \sum_{k \in K}\sum_{s\in S} w_{r_i, r_j}^{s,k} \Biggr)}{\ \forall r_i \in \R, r_j \in \precedenceSet{r_i}^{o,o}}\label{model:subline:pax:time:3}\\ 
        &{b_{r_i} + \dropoff{r_i} - \pickup{r_j}}{\leq M_1 \cdot \Biggl( 1 - \sum_{k \in K}\sum_{s\in S} w_{r_i, r_j}^{s,k} \Biggr)}{\ \forall r_i \in \R, r_j \in \precedenceSet{r_i}^{o,d}}\label{model:subline:pax:time:4}\\ 
        &{b_{r_i} + \dropoff{r_i} - \dropoff{r_j}}{\leq M_1 \cdot \Biggl( 1 - \sum_{k \in K}\sum_{s\in S} w_{r_i, r_j}^{s,k} \Biggr)}{\ \forall r_i \in \R, r_j \in \precedenceSet{r_i}^{d,d}}\label{model:subline:pax:time:5}\\ 
        &{\depS_{o_r}^{s,k}}{\geq \assignS_r^{s,k} \cdot ( \pickup{r} + b_r)}{\qquad\forall\ \allR{r}}\label{model:subline:pax:time:6}\\
        &{\depS_{d_r}^{s,k}}{\geq \assignS_r^{s,k} \cdot ( \dropoff{r} + b_r )}{\qquad\forall\ \allR{r}}\label{model:subline:pax:time:7}\\
        &{ e_{o_r}}{\leq \pickup{r} \leq l_{o_r}}{\qquad\forall\ \allR{r}}\label{model:subline:pax:time:8}\\ 
        &{e_{d_r}}{\leq \dropoff{r} \leq l_{d_r}}{\qquad\forall\ \allR{r}}\label{model:subline:pax:time:9}\\
        &{\dropoff{r} + b_r - \pickup{r}}{ \leq \alpha \cdot t_{\org{r}, \dest{r}}}{\qquad \forall\ \allR{r}}\label{model:subline:pax:time:10} 
    \end{align}
\end{subequations}

Constraint \eqref{model:subline:subline:1} ensures each vehicle is used at most once and constraint \eqref{model:subline:subline:2} ensures that a started vehicle also ends at a bus stop. Constraint \eqref{model:subline:subline:3} makes sure every started vehicle turns around at some station, while \eqref{model:subline:subline:4} denotes that the last subline of each vehicle does not turn. Constraints \eqref{model:subline:subline:5} to \eqref{model:subline:subline:9} are for flow conservation between sublines and their turn stops. The last subline ends at its end stop, which is controlled by \eqref{model:subline:subline:10}. Constraint \eqref{model:subline:subline:11} tracks the amount of vehicles which are used. The upper capacity of every vehicle is ensured by \eqref{model:subline:subline:12} and \eqref{model:subline:subline:13}, for both subline directions. Constraints \eqref{model:subline:subline:14} to \eqref{model:subline:subline:18} ensure that only vehicles which are started can travel to and between stops. Constraint \eqref{model:subline:subline:19} is a symmetry breaking constraint which says that vehicles with a smaller index start at a smaller station.

Constraint \eqref{model:subline:time:1} ensures a vehicle departs from a bus stop only after it has arrived, while constraints \eqref{model:subline:time:2} and \eqref{model:subline:time:3} ensure the arrival time at the next bus stop respects the minimum travel time. Constraints \eqref{model:subline:time:4} and \eqref{model:subline:time:5} take into account the turning time, linking a subline's end time with the subsequent subline's start time at the same stop. Finally, constraints \eqref{model:subline:time:6} ensures all times are positive. 

Constraint \eqref{model:subline:pax:subline:1} ensures a passenger is picked up at most once and, by constraint \eqref{model:subline:pax:subline:2}, only if the subline they are assigned to also stops at their origin and destination. The link between earliest pick-up and departure from the origin times, as well as latest drop-off and arrival at the destination times, is handled with constraints \eqref{model:subline:pax:subline:5} and \eqref{model:subline:pax:subline:6}, respectively. Finally, constraints \eqref{model:subline:pax:subline:7} and \eqref{model:subline:pax:subline:8} link the variable $w_{r_i, r_j}$ to denote if two passengers are assigned to the same subline of the same vehicle.

Constraints \eqref{model:subline:pax:time:1} and \eqref{model:subline:pax:time:2} place a lower bound on the pick-up time and drop-off time of each passenger, dependent on the vehicle's arrival time at the corresponding station. Constraints \eqref{model:subline:pax:time:3} to \eqref{model:subline:pax:time:5} ensure the precedence rules for boarding are respected and add sufficient service times between serving customers. Then, constraints \eqref{model:subline:pax:time:6} and \eqref{model:subline:pax:time:7} ensures the vehicle can only depart after the last passenger has fully boarded or alighted. The time windows of each passenger is guaranteed by \eqref{model:subline:pax:time:8} and \eqref{model:subline:pax:time:9}, while the maximum travel time is limited by constraint \eqref{model:subline:pax:time:10}.

\subsection{Location-Based Formulation}\label{appendix:milp:darp}
In this section, we describe the construction of the underlying graph for the Location-Based formulation in more detail, as well as presenting the full MILP model. All notation is summarized in \Cref{tab:DARP:Vars}.

\begin{table}[ht!]
    \caption{Summary of notation for the Location-Based model.}
    \label{tab:DARP:Vars}
    \begin{tabularx}{\textwidth}{ll}
    \toprule
    Notation & Definition\\
    \midrule
    \multicolumn{2}{l}{Parameters} \\
    \midrule
    $\dstart$& start depot  \\
    $\dend$ & end depot  \\
    $P$ & set of pick-up nodes\\
    $D$ & set of delivery node\\
    $\barP$ & set of start-turn nodes, before a pick-up node\\
    $\barD$ & set of end-turn nodes, after a drop-off node\\
    $\NR$ & set of all pick-up and delivery nodes in all directions, depending on requests $\R$ \\
    $\HR$ & set of all pick-up, delivery, and depot nodes, depending on requests $\R$ \\
    $\ER$ & set of all edges between nodes in $\ER$ \\
    \midrule
    \multicolumn{2}{l}{Binary Decision Variables} \\
    \midrule
    $x_{i,j}^k$  & $1$ if vehicle $k$ travels on arc $(i,j) \in \ER$\\
    $z_k$ & $1$ if vehicle $k$ is used \\
    \midrule
    \multicolumn{2}{l}{Continuous Decision Variables} \\
    \midrule
    $B_i$ & start of service time at bus stop $i$\\
    $Q_i$ & passenger load departing bus stop $i$ \\
    $L_r$ & ride time of passenger $r$ \\
       \bottomrule
    \end{tabularx}
\end{table}

We introduce \emph{ascending} bus stops $\Hleft := \{h_1, \ldots, h_n\}$ and \emph{descending} bus stops $\Hright := \{h_{n+1}, \ldots, h_{2n}\}$. Bus stops $h_i \in \Hright$ and $h_{n+i}\in \Hleft$ are virtual copies of stop $i \in H$. In reality, these may be the same stop on opposite sides of the road, for vehicles travelling in opposite directions. We connect all bus stops $h_i \in \Hright$ with their corresponding $h_{n+i} \in \Hleft$, in both directions (corresponding to a turn at station $i$), and enforce that two stops in either set can only be served in upstream order with respect to the line.

Similar to the classical DARP formulation, for each request $r \in \R$, we construct four nodes $o_r, d_r, \baro_r, \bard_r$ in a liDARP-Graph $\GR = (\HR, \ER)$. Here, the nodes $o_r$ and $d_r$ correspond to the classic pick-up and drop-off nodes of $r$. The node $\baro_r$ is a \emph{start-turn} node, denoting that the vehicle is at $o_r$, facing the opposite direction, i.e., before it turns and picks up request $r$ at $o_r$. Similarly, the node $\bard_r$ is an \emph{end-turn} node, denoting that the vehicle is at $d_r$, has dropped off the request $r$, and is now turning to continue in the opposite direction. Then,
\begin{itemize}
    \item if $r \in \Rleft$: we construct $o_r$ at the origin $h_{o_r} \in \Hleft$, $d_r$ at the destination $h_{d_r} \in \Hleft$ stop, $\baro_r$ at $h_{n + o_r} \in \Hright$, and $\bard_r$ at $h_{n + d_r} \in \Hright$.
    \item if $r \in \Rright$: we construct $o_r$ at the origin $h_{o_r} \in \Hright$, $d_r$ at the destination $h_{d_r} \in \Hright$ stop, $\baro_r$ at $h_{o_r - n} \in \Hleft$, and $\bard_i$ at $h_{d_r - n} \in \Hleft$.
\end{itemize}

We set $P := \{o_r : r \in \R \}$,  $D := \{d_r: r \in \R\}$,  $\barP := \{\baro_r :r \in \R \}$, and  $\barD := \{\bard_r: r \in \R\}$ to denote the sets of these bus stops. We additionally define time windows on the nodes $\baro_r$ and $\bard_r$ of every request $r \in \R$, dependent on the time windows on $o_r$ and $d_r$, respectively, accounting for the boarding and turn times.

We introduce two depots, the \emph{start depot} $\dstart$ and the \emph{end depot} $\dend$, where vehicles start and end their route. Let $\NR := P \cup D \cup \barP \cup \barD$ denote all pick-up, drop-off, and turn stops, and $\HR :=  \{\dstart, \dend\} \cup \NR$ all nodes including the depots.

The edge set $\ER:= \bigcup_{i = 1}^{10} \ER^i$ between nodes in $\HR$ is constructed as follows:
\begin{itemize}
    \item $\ER^1 := \{(o_r, d_r) \in P \times D : r \in \R \}$, connecting each request's origin with its destination,
    \item $\ER^2 := \{(v_i,w_j) \in (P \cup D \cup \barD) \times (P \cup D \cup \barP): h_{v_i}, h_{w_j} \in \Hleft, i \neq j, v_i \text{ precedes } w_j, i,j \in \R\}$, connecting all \emph{ascending} stops to subsequent stops in the same direction,
    \item $\ER^3 := \{(v_i,w_j) \in (P \cup D \cup \barD) \times (P \cup D \cup \barP): h_{v_i}, h_{w_j} \in \Hright, i \neq j, v_i \text{ precedes } w_j, i,j \in \R\}$, connecting all \emph{descending}  stops to subsequent stops in the same direction,
    \item $\ER^4:= \{(\baro_r, o_r) \in \barP \times P : r \in \R \}$, connecting the start-turn stop of each request with its corresponding origin stop in the opposite direction,
    \item $\ER^5:= \{(d_r, \bard_r) \in D \times \barD: r \in \R \}$, connecting the destination stop of each request with its corresponding end-turn stop in the opposite direction,
    \item $\ER^6:= \{(v_i,w_j) \in (P \cup D \cup \barD) \times (P \cup D \cup \barP) : v_i,w_j \in \Hleft, h_{v_i} = h_{w_j},i \neq j,\allowbreak
    \lnot (v_i \in D \land w_j \in P), e_{v_i} \leq e_{w_j}, l_{v_i} \leq l_{w_j}, i,j \in \R \}$, connecting \emph{ascending} stops at the same original physical bus stop if they are compatible regarding their time windows,
    \item $\ER^7:= \{(v_i,w_j) \in (P \cup D \cup \barD) \times (P \cup D \cup \barP): v_i,w_j \in \Hright, h_{v_i} = h_{w_j},i \neq j, \allowbreak
    \lnot (v_i \in P \land w_j \in D),  e_{v_i} \leq e_{w_j}, l_{v_i} \leq l_{w_j}, i,j \in \R \}$, connecting \emph{descending} stops at the same original physical bus stop if they are compatible regarding their time windows,
    \item $\ER^8 := \{ (\dstart, \dend)\}$, connecting the starting depot to the ending depot to allow for unused vehicles,
    \item $\ER^9 := \{(\dstart, o_r) \in \{\dstart\} \times P: r \in \R\} \cup \{(\dstart, \baro_r) \in \{\dstart\} \times \barP: r \in \R\}$, connecting the start depot to all pick-up locations and their start-turn stops,
    \item $\ER^{10} := \{(d_r, \dend) \in D \times \{\dend\} : r \in \R \} \cup \{(\bard_r, \dend) \in \barD \times \{\dend\} : r \in \R \}$, connecting all drop-off locations and their end-turn stops to the end depot.
\end{itemize}

Here, we write \emph{$v$ precedes $w$} to denote that the bus stations corresponding to $v$ precedes that corresponding to $w$ with respect to the corresponding line direction. We use $\Eturn := \ER^4 \cup \ER^5$ to denote all the edges on which the vehicles turn. Each edge is only added once, even if it appears in multiple sets.
The travel time of the edges is given by the original network, where a turn takes $\tturn$ and travel between two pick-ups or drop-offs at the same physical stop is instantaneous.

In our model, the binary variable $x_{i,j}^k$ denotes if a vehicle $k$ travels on defined arcs $(i,j) \in \ER$. Variables $z_k$ denote if vehicle $k$ is used in the solution.

We define that the vehicle loads $q_{\dstart} := q_{\dend} = 0$, $q_i := 1$ for all $i \in P$ and $q_i := -1$ for all $i \in D$. Additionally, we set the service times $b_{\dstart} := b_{\dend} = 0$ and $b_i = 0$ for all $i \in \barP \cup \barD$.
Let $Q_i$ denote the passenger load of a vehicle departing a stop $i$ and let the continuous variable $B_i$ denote the start of service time at stop $i$. Note that these do not require an index for the vehicle $k$ as each node can be visited by at most one vehicle and the vehicles are homogeneous with a maximum capacity $\Qmax$, as has been discussed in \cite{cordeau_branch-and-cut_2006}. We require that $Q_{\dstart}^k := Q_{\dend}^k = 0$ for all $k \in K$, thus removing these variables from the model. 

To strengthen the model, we introduce a symmetry breaking constraint which enforces that vehicles of lower index are used first. Then, the full Location-Based model is given by:
\begin{maxi!}[2]
    {x}{ c_1 \Biggl(\sum_{k \in K} \sum_{i \in P} t_{i, i+m} \cdot x_{i, i+m}^k - \sum_{k \in K} \sum_{\subjH{i}{j}} \sum_{i \in \HR} t_{i,j} \cdot x_{i,j}^k \Biggr)  
    + c_2 \sum_{k \in K} \sum_{\subjH{i}{j}} \sum_{i \in P} x_{i,j}^k }{}{\notag}
    \addConstraint{\sum_{k \in K} \sum_{\subjH{i}{j}} x_{i,j}^k}{\leq 1 \qquad\forall\ \allP{i}}{\label{model:darp:1}} 
    \addConstraint{\sum_{\subjH{i}{j}} x_{i,j}^k - \sum_{\subjH{m+i}{j}} x_{m{+}i,j}^k}{= 0 \qquad\forall\ \allP{i}, \allK}{\label{model:darp:2}} 
    \addConstraint{\sum_{j \in P \cup \barP} x_{\dstart, j}^k}{
    = 1 \qquad\forall\ \allK}{\label{model:darp:3}} 
    \addConstraint{\sum_{i \in D \cup \barD} x_{i, \dend}^k}{= 1 \qquad\forall\ \allK}{\label{model:darp:4}} 
    \addConstraint{\sum_{\subjH{j}{i}} x_{j,i}^k - \sum_{\subjH{i}{j}}x_{i,j}^k}{= 0 \qquad\forall\ \allNR{i}, \allK}{\label{model:darp:5}} 
    \addConstraint{1 - \sum_{j\in \NR} \sum_{\substack{i \in \NR: \\ (i,j) \in \ER}} x_{i,j}^k}{\leq M_3\cdot x_{\dstart, \dend}^k \qquad\forall\ \allK}{\label{model:darp:6}} 
    \addConstraint{B_j}{\geq (B_{\dstart}^k + b_{\dstart} + t_{\dstart,j}) \cdot x_{\dstart, j}^k \qquad \forall\ j \in P \cup \barP \cup \{\dend\}, \allK}{\label{model:darp:7}} 
    \addConstraint{B_{\dend}^k}{\geq (B_i + b_i + t_{i, \dend}) \cdot x_{i, \dend}^k \qquad \forall\ i \in D \cup \barD \cup \{\dstart\}, \allK}{\label{model:darp:8}}
    \addConstraint{B_j}{\geq (B_i + b_i + t_{i,j}) \cdot \sum_{k \in K} x_{i,j}^k \qquad\forall\ \allNR{i,j} : (i,j) \in \ER}{\label{model:darp:9}}
    \addConstraint{B_i}{\geq e_i + \sum_{j \in \HR \setminus \{i\}} \Biggl( \max\{0, e_j - e_i + b_j + t_{j,i}\} \cdot \sum_{k \in K} x_{j,i}^k \Biggr) \qquad\forall\ \allNR{i}}{\label{model:darp:10}} 
    \addConstraint{B_i}{\leq l_i - \sum_{j \in \HR \setminus \{i\}} \Biggl( \max \{0, l_i - l_j + b_i + t_{i,j}\} \cdot \sum_{k \in K} x_{i,j}^k \Biggr) \qquad\forall\ \allNR{i}}{\label{model:darp:11}}
    \addConstraint{L_i}{= B_{i+m} - (B_i + b_i) \qquad\forall\ \allP{i}}{\label{model:darp:12}} 
    \addConstraint{t_{i, i{+}m}}{\leq L_i \leq \alpha \cdot t_{i, i{+}m} \qquad\forall\ \allP{i}}{\label{model:darp:13}} 
    \addConstraint{Q_j}{\geq q_j \cdot x_{\dstart,j}^k \qquad \forall\ j \in P \cup \barP, \allK}{\label{model:darp:14}} 
    \addConstraint{Q_j}{\geq (Q_i + q_j) \cdot \sum_{k \in K} x_{i,j}^k \qquad\forall\ \allNR{i,j} : (i,j) \in \ER}{\label{model:darp:15}}
    \addConstraint{0}{\geq Q_i \cdot x_{i, \dend}^k \qquad \forall\ i \in D \cup \barD, \allK}{\label{model:darp:16}}
    \addConstraint{Q_i}{\leq \Qmax \cdot (1 - \sum_{k \in K} x_{i,j}^k) \qquad \forall\ (i,j) \in \Eturn}{\label{model:darp:17}} 
    \addConstraint{Q_i}{\geq -\Qmax \cdot (1 - \sum_{k \in K} x_{i,j}) \qquad \forall\ (i,j) \in \Eturn}{\label{model:darp:18}}
    \addConstraint{Q_i}{\leq \Qmax \cdot \sum_{k \in K} \sum_{\subjH{i}{j}} x_{i,j}^k \qquad\forall\ \allHR{i}}{\label{model:darp:19}} 
    \addConstraint{z_k}{\geq \frac{1}{|\HR|^2} \cdot \sum_{(i,j) \in \ER} x_{i,j}^k \qquad\forall\ \allK}{\label{model:darp:20}} 
    \addConstraint{z_{k}}{\geq z_{k+1} \qquad\forall\ \allK}{\label{model:darp:21}} 
\end{maxi!}

Constraints \eqref{model:darp:1} and \eqref{model:darp:2} ensure each passenger is picked up at most once and is dropped-off by the same vehicle. Vehicles must start \eqref{model:darp:3} at $\dstart$ and end at $\dend$ \eqref{model:darp:4} depots, maintaining flow conservation across all arcs \eqref{model:darp:5}. Only unused vehicles may use the arc $(\dstart, \dend)$, as denoted by \eqref{model:darp:6}, where we use a big-M constraint with $M_3 := |\ER|$. The service start times for leaving and entering the depot, as well as consistency across arcs, is handled by \eqref{model:darp:7}--\eqref{model:darp:9}. Constraints \eqref{model:darp:10} and \eqref{model:darp:11} ensure time consistency regarding the requests time windows. The maximum ride time of each request is defined and bounded by constraints \eqref{model:darp:12} and \eqref{model:darp:13}. Load constraints \eqref{model:darp:14}--\eqref{model:darp:19} ensure vehicles respect capacity limits at each bus stop as well as on turning arcs. Constraint \eqref{model:darp:20} counts the number of required vehicles. Finally, we use the symmetry breaking constraint \eqref{model:darp:21} to improve computational times.
\end{document}